# Capacitated Bounded Cardinality Hub Routing Problem: Model and Solution Algorithm

Shahin Gelareh[a,b], Rahimeh Neamatian Monemi[c], Frédéric Semet[d,*]

[a] *LGI2A (EA 3926), Université d'Artois, F-62400 Béthune, France*
[b] *Univ Lille Nord de France, F-59000 Lille, France*
[c] *Ifsttar, Univ. Lille Nord de France, rue Elísée Reclus 20, 59666 Villeneuve d'Ascq, France*
[d] *CRIStAL, UMR 9189-CNRS, Ecole Centrale de Lille, 59651 Villeneuve d'Ascq, France*

**Abstract**

In this paper, we address the Bounded Cardinality Hub Location Routing with Route Capacity wherein each hub acts as a transshipment node for one directed route. The number of hubs lies between a minimum and a maximum and the hub-level network is a complete subgraph. The transshipment operations take place at the hub nodes and flow transfer time from a hub-level transporter to a spoke-level vehicle influences spoke-to-hub allocations. We propose a mathematical model and a branch-and-cut algorithm based on Benders decomposition to solve the problem. To accelerate convergence, our solution framework embeds an efficient heuristic producing high-quality solutions in short computation times. In addition, we show how symmetry can be exploited to accelerate and improve the performance of our method.

*Keywords:* Hub Location, Location Routing, Branch-and-Bound, Benders Decomposition

## 1. Introduction

Many transportation operators and logistics service providers operate on hub-and-spoke structures. In this way, the long haul transport of their cargo at national or international scale is done by using larger transporters –with higher volume concentration– circulating between major hubs. At a short-distance level, transportation service is realized by means of smaller vehicles. In such services, the long-haul transport services are in fact some point-to-point shuttle services calling at both end-point hub nodes as long as this is feasible with respect to the geographical and environmental restrictions. On the other hand, at a hub, the spoke node customers allocated to the hub are served on service routes with regular, reliable and often fixed frequency (for example in liner shipping, barge services etc.). On such routes, a feeder service (spoke-level route) is in charge of distributing the import volumes destined to all the spoke nodes and at the same time collecting the export volumes from the spoke nodes to be delivered to the hub node for further distribution operations. Often such feeder services are directed and there are no more than one route along which the spoke nodes of a hub are served. In addition, the intersection of two feeder routes associated to two distinct hub nodes is often an empty set due to geographical distribution of nodes.

In this paper, we address precisely a such problem called the Bounded Cardinality Capacitated Hub Routing Problem (BCCHRP). The BCCHRP is defined on a time-weighted graph $G(V, A, T)$, where $V = \{1, 2, \ldots, n\}$ is the set of nodes, $A$ is the set of arcs and $T = [t_a : a \in A]$ is the set of weights. $\alpha$ is the mode-dependent factor of economies of scale (hub-level efficiency) applicable to the vehicles traversing the hub-level edges, $\Gamma$ represents the minimum number of customers along every feeder route and $C$ is the capacity of vehicles on the feeder routes. We seek a partition of $V$ into a set of hub nodes $H$ and a set of spoke nodes associated with each hub node, $S^h : h \in H$, minimizing the total travel time. The travel time for a given origin-destination

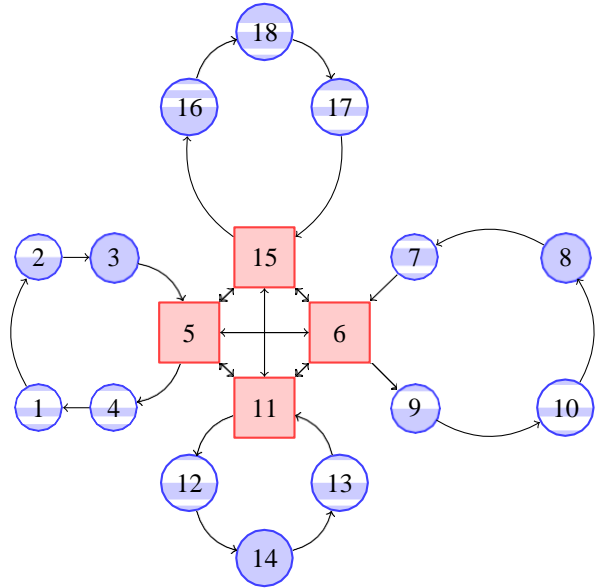

**Figure 1:** An example of the network structure of BCCHRP. Rectangles represent hub nodes and circles represent spoke nodes. The routes connecting sets of spoke nodes correspond to circuits. The hub-level network is a complete bidirectional graph.

(O-D) flow accounts for the travel time from origin $i$ on the feeder route up until the hub node $k$ where $i \in S^k$ (if $i$ is not a hub itself), plus the transshipment time at the hub port $k$, plus the travel time $\alpha t_{k,l}$ (discounted travel time due to the use of faster lines) on the direct edge between hub $k$ and hub $l$, plus the transshipment time at the hub node $l$ where $j \in S^l$, plus the travel time from the hub node $l$ to the destination $j$. Note that for an origin or destination, which is a hub node, neither transshipment time nor feeder-level transportation times are taken into account when calculating the objective function for the hub end-node(s).

Figure 1 depicts the structure of such networks. In this figure, $|H| = |\{5,6,11,15\}| = 4$, $S^5 = \{1,2,3,4\}$, $S^6 = \{7,8,9,10\}$, $S^{11} = \{12,13,14\}$, $S^{15} = \{16,17,18\}$ and $|V| = 18$. The flows over routes associated to each rectangular hub node respect vehicle capacity. Note that the routes are circuits.

Hub Location Problems (HLPs) originate from the seminal work of Hakimi (1964) on finding the optimal location of a single switching center that minimizes the total wire length in a communication network. He showed that one can limit oneself to find the vertex median of the corresponding graph. Later, Goldman (1969) proved that Hakimi's argument holds even for more general cases. While studies on the hub-and-spoke network structures in transportation dates back to Toh and Higgins (1985), the models in O'Kelly (1986a,b) were the first continuous location models for the HLPs dealing with locating two hub nodes in the plane. Later, starting from a quadratic integer programming model in O'Kelly (1987), the research efforts were mostly focused on discrete variants.

In classical models, the following assumptions were almost always considered: (1) The hub level network was a complete graph, (2) There was a discount factor associated with the use of hub edges reflecting economies of scale, (3) Direct connections between non-hub nodes were not allowed, (4) Costs were proportional to the distances —i.e. the triangle inequality held, and (5) All nodes were candidates to become hubs. As a result, (1), (3) and (4) together were implying that every origin-destination flow passes through one or at most two hubs. However, to the best of our knowledge, for the first time Nickel et al. (2001) and later in Yoon and Current (2008), Gelareh (2008); Gelareh and Nickel (2011), Contreras et al. (2009, 2010) and Alumur and Kara (2009); Alumur et al. (2009) proposed models that relax these assumptions and proposed better approximation of real practice in different applications in freight transportation and telecommunications. There are also other contributions that deal with the flow-dependent costs —for example (O'Kelly and Bryan, 1998) and/or congestion (de Camargo et al., 2009a), among others.

Interested readers are also referred to the recent HLPs literature review in Alumur and Kara (2008),

Campbell and O'Kelly (2012) and Farahani et al. (2013) as well as other contributions by Campbell et al. (2002) and Kara and Taner (2011).

In this paper, we deal with a hub-and-spoke network structure wherein the hub-level network is complete, the allocation of spoke nodes to the hub nodes follows a single allocation scheme, and more importantly the spoke nodes allocated to a given hub node are served on a directed route. The problem can be classified as a variant of *Hub Location and Routing Problems (HLRPs)*.

Gelareh et al. (2013) proposed a mixed integer linear programming formulation for simultaneous network design and fleet deployment of a deep-sea liner service provider. The underlying network design problem is based on a 4-index formulation of the hub location problem (mainly models in Gelareh and Nickel (2011) and Gelareh et al. (2010)) and the spoke nodes allocated to a hub node form directed cycles (strings). The vessels class and quantity on every route needs to be determined from a given fleet of vessels with different capacities and different factors of economies of scale. The objective is to minimize a weighted sum of transit times, and fixed deployment costs. While the existing general-purpose MIP solvers were unable to solve even very small size problem instances in a reasonable time, the authors proposed a Lagrangian decomposition approach equipped with a heuristic procedure that shown to be efficient. They applied their model to a case study from the liner shipping industry. Gelareh et al. (2013) referred to this problem as *p-String Planning Problem (pSPP)*.

Other related contributions include Çetiner et al. (2010) for a multiple allocation hub location and routing model applied to postal delivery applications. Two objective functions were taken into account and the vehicles performing routes had an infinite capacity. Again, for the similar applications, Wasner and Zäpfel (2004) proposed a model wherein direct connections between spoke nodes were allowed. de Camargo et al. (2013) proposed a similar model when the route length is bounded and a Benders decomposition approach for its solution. Nagy and Salhi (1998) proposed a hub location routing with capacity constraints. In their model, each spoke nodes may be visited twice since pick-up and delivery services are performed on two distinct routes. Rodríguez-Martín et al. (2014) referred to this class of problems by *Hub Location and Routing*. They proposed a mathematical model and a branch-and-cut approach to a solve some instances of up to 50 nodes. The model and the cuts are based on a previous work on the Plant Cycle Location Problem (PCLP) Labbé et al. (2004). In this work, the number of hub nodes is fixed and the capacity constraint is on the maximum number of spoke nodes served on a route.

The current work differs from Rodríguez-Martín et al. (2014) by the following assumptions: 1) there is a minimum and a maximum number of hub nodes that can be installed while in Rodríguez-Martín et al. (2014) this is a fixed number, 2) in this work our objective is to minimize the total transportation and transshipment times while in Rodríguez-Martín et al. (2014), the objective is to only minimize transportation cost, 3) we consider directed routes while in Rodríguez-Martín et al. (2014) the routes are undirected and, 4) here, capacity constraints are imposed while in Rodríguez-Martín et al. (2014) the number of spoke nodes along a route is bounded. Table 1 summarizes some key features in closely related contributions.

**Table 1:** A summary of main elements of the relevant contributions in literature.

| Work | Allocation Scheme | No. Hubs | Objective | Capacity | Cycle Length | No. Vehicles | Solution method |
|---|---|---|---|---|---|---|---|
| Nagy and Salhi (1998) | pickup/delivery routes | endogenous | cost | yes | yes | endogenous | MIP + heuristic |
| Cetiner et al. (2006) | multiple allocation | endogenous | cost+fleet | no | yes | endogenous | heuristic |
| de Camargo et al. (2013) | single allocation | endogenous | cost | no | yes | endogenous | MIP + Benders decomposition |
| Wasner and Zäpfel (2004) | multiple allocation | endogenous | cost | yes | yes | endogenous | MIP + heuristic |
| Rodríguez-Martín et al. (2014) | single allocation | exogenous | cost | #. spoke per route | no | one per route | MIP + branch-and-cut |
| Gelareh et al. (2013) | single allocation | exogenous | cost+fleet | yes | multiple of weeks | variable | MIP+Lagrangian decomposition |
| This work | single allocation | $q = 3 \leq .. \leq p$ | time (transit+transshipment) | yes | $\geq 2$ spokes | one per route | MIP+branch-and-cut+Benders |

The main contributions of this work are threefold. First, we propose a mixed integer linear formulation with 2-index design variables. Then we describe an efficient exact branch-and-cut method based on a Benders decomposition for this problem. Last we propose a systematic way to exploit symmetry in generating (relative) interior points of the polytope to identify sharper Benders cuts.

The remainder of this paper is organized as follows. In section 2, the mixed integer linear programming model is given. In section 3, the branch-and-cut method based on a Benders decomposition is presented. Computational experiments and discussions are reported in section 4. Finally, in section 5, we summarize, draw conclusions and provide suggestions for further research directions.

## 2. Mathematical model

We propose a mixed integer linear models with 2-index design variables. The parameters used in our model are listed in Table 2:

**Table 2:** Model Parameters.

| | |
|---|---|
| $w_{ij}$: | the volume of goods transferred from $i$ to $j$, |
| $t_{ij}$: | the time associated with the arc $(i,j)$, |
| $\alpha$: | the factor of economies of scale (the factor of travel time efficiency over hub arcs), |
| $p$: | the upper bound on the number of hub nodes, |
| $q$: | the lower bound on the number of hub nodes, |
| $\Gamma$: | the minimum number of spoke nodes allocated to each hub node, |
| $C$: | the vehicle capacity, |
| $\varphi^k$: | the (fixed) average transshipment time at hub $k$. |

We assume that the transshipment time $\varphi^k$ at hub $k$ is independent of the type of goods and of the volumes.

The decision variables are the followings: $x_{ijkl}$ is the fraction of flow from $i$ to $j$ traversing inter-hub edge $\{k,l\}$. $s_{ijkl}$ represents the fraction of flow from $i$ to $j$ traversing non-hub edge $\{k,l\}$. $r_{ij}$ is equal to 1, if the arc $(i,j)$ belongs to a route, 0 otherwise. $z_{ik}$ is equal to 1, if the spoke node $i$ is allocated to the hub node $k$, 0 otherwise.

### 2.1. A compact formulation for the BCCHRP

The BCCHRP can be modeled as follows:

$$(BCCHRP) \min \sum_{i,j,k,l} t_{kl}(s_{ijkl} + \alpha x_{ijkl}) + \sum_{i,j,k,l:(k \neq i \vee l \neq j)} (\varphi^k + \varphi^l) x_{ijkl} \tag{1}$$

$$s.t.$$

$$q \leq \sum_k z_{kk} \leq p \tag{2}$$

$$\sum_l z_{kl} = 1 \qquad \forall k \in V \tag{3}$$

$$z_{ik} \leq z_{kk} \qquad \forall i,k \in V : k \neq i \tag{4}$$

$$\sum_i z_{ik} \geq \Gamma z_{kk} \qquad \forall k \in V \tag{5}$$

$$\sum_{j \neq i} r_{ij} = 1 \qquad \forall i \in V \qquad (6)$$

$$\sum_{j \neq i} r_{ji} = 1 \qquad \forall i \in V \qquad (7)$$

$$r_{ij} + r_{ji} \leq 2 - z_{ik} - z_{jl} \qquad \forall i, j, k, l \in V : j \neq i, k \neq l \qquad (8)$$

$$r_{ij} + r_{ji} \leq 1 \qquad \forall i, j \in V : j \neq i, \qquad (9)$$

$$\sum_{k \neq i}(x_{ijik} + s_{ijik}) = 1, \qquad \forall i, j \in V : j \neq i, \qquad (10)$$

$$\sum_{l \neq j}(x_{ijlj} + s_{ijlj}) = 1, \qquad \forall i, j \in V : j \neq i, \qquad (11)$$

$$\sum_{l \neq i,k}(x_{ijkl} + s_{ijkl}) = \sum_{l \neq j,k}(x_{ijlk} + s_{ijlk}), \qquad \forall i, j, k \in V, k \notin \{i, j\}, \qquad (12)$$

$$\sum_{l \neq k} x_{ijkl} \leq z_{kk} \qquad \forall i, j, k \in V : j \neq i, k < l \qquad (13)$$

$$\sum_{l \neq k} x_{ijlk} \leq z_{kk} \qquad \forall i, j, k \in V : j \neq i, k < l \qquad (14)$$

$$s_{ijkl} \leq r_{kl} \qquad \forall i, j, k, l \in V : l \neq k \qquad (15)$$

$$\sum_{ijl: j \neq i} w_{ij} s_{ijkl} \leq C, \qquad \forall k \in V \qquad (16)$$

$$r \in \mathbb{B}^{|V|^2}, z \in \mathbb{B}^{|V| \times |V|}, x_{ijkl}, s_{ijkl} \in \mathbb{R}_{[0,1]}^{|V|^4} \qquad (17)$$

The objective function (1) accounts for the total transportation and transshipment times. The transportation part is comprised of travel times on the spoke links plus the hub-level edges, which provide faster services (α-discounted times). The transshipment time is considered for a flow that traverses hub edges where at least one of the endpoints of this hub edge is different from the origin or destination of the flow. No transshipment costs is applicable at an origin (destination) that is a hub node by itself. Constraint (2) sets bounds on the number of hub nodes that can be installed while constraints (3) guarantee that every node is allocated to exactly one hub. A self-allocation at node $i$ (i.e. $z_{ii} = 1$) indicates that node $i$ is a hub node. A spoke node can only be allocated to a hub node as stated in constraints (4). If a node is designated as a hub node, there must be at least $\Gamma$ nodes (including itself) that are allocated to it, as stated in constraint (5). Constraints (6) and (7) ensure that one non-hub arc leaves and one non-hub arc arrives at every node, respectively. Constraints (8) ensure that a spoke arc cannot exist if its endpoints are allocated to different hubs and constraints (9) guarantee that only one arc between two nodes is in the solution. This constraint guarantees the feasibility of the solution with respect to the pick-up and delivery services.

Constraints (10)-(14) are related to the flow conservation. Constraints (10)-(11) are the out-going and in-going flow restrictions. Constraints (12) are the flow conservation restrictions. Constraints (13)-(14) impose that the flow from node $i$ to node $j$ is routed through node $k$ if and only if it is selected as a hub node. Constraints (15) make sure that a spoke flow will traverse an existing spoke arc. Constraints (16) are the capacity restrictions.

### 2.2. Valid inequalities

In this section, we present valid inequalities to strengthen the linear relaxation of (BCCHRP). We note that $r(\delta(S)) = \sum_{i \in S, j \in V/S} r_{ij}$, $r(A(S)) = \sum_{(i,j) \in A(S)} r_{ij}$ where $A(S) = \{a = (i, j) \in A : i, j \in S\}$ and $r(P) = \sum_{(i,j) \in P} r_{ij}$ where $P = \{(i_0, i_1), (i_1, i_2), \ldots, (i_{k-1}, i_k)\}$.

**Proposition 1.** *If $i \in S \subset V$ and $i$ is a spoke node allocated to a hub node in $V/S$, there is at least one arc going out of (entering into) $S$:*

$$a) \quad r(\delta^+(S)) \geq \sum_{j \in V/S} z_{ij}, \qquad \forall i \in S \subset V \qquad (18)$$

$$b) \ r(\delta^-(S)) \geq \sum_{j \in V/S} z_{ij}, \qquad \forall i \in S \subset V. \qquad (19)$$

*Proof.* a) For a given $S \subset V$ and any $i \in S$, if no arc $(i, j)$ with $j \in V/S$ exists, therefore there exists a path from $i$ to the hub node where $i$ is allocated to it such that this path does not cut $S$ and remains entirely inside $S$. Thus, there exists no $j$ for which $z_{ij} = 1$ for any choice of $j \in V/S$. On the other hand, if $i$ is allocated to a hub node in $V/S$ there is a unique path cutting $S$ (at least once) towards the hub node where $i$ is allocated to it.

b) Analogous to the proof of a). □

Rodríguez-Martín et al. (2014) and Labbé et al. (2004) propose a used a similar set of valid inequalities for the case of undirected graphs.

**Proposition 2.** *If $i \in S \subset V$ and $i$ is a spoke node allocated to a hub node in $V/S$ and $i' \in V/S$ is a spoke node allocated to a hub node in $S$, we have:*

$$r(\delta^+(S)) \geq \sum_{j \in V/S} z_{ij} + \sum_{j' \in S} z_{i'j'}, \qquad \forall i \in S \subset V, i' \in V/S \qquad (20)$$

$$r(\delta^-(S)) \geq \sum_{j \in V/S} z_{ij} + \sum_{j' \in S} z_{i'j'}, \qquad \forall i \in S \subset V, i' \in V/S \qquad (21)$$

*Proof.* The proof is very similar to Proposition 1. Rodríguez-Martín et al. (2014) and Labbé et al. (2004) proposed a similar class of valid inequalities for their problems. □

**Proposition 3.** *If $(i - j)(or(j, i))$ is an arc and $j$ allocated to a hub $k$ then $i$ must be also allocated to the same hub $k$:*

$$r_{ij} + r_{ji} + z_{jk} \leq 1 + z_{ik} \qquad \forall i, j, k : j \varsubsetneq i \qquad (22)$$

**Proposition 4.** *For all $S \subseteq V$,*

$$r(A(S)) - \sum_{i \in S} z_{ii} \leq |S| - \lceil \frac{|S|}{U} \rceil \qquad (23)$$

*is a valid inequality.*

*Proof.* Given that $r(E(S)) = r(A(S))$, c.f. Proposition 2 in Rodríguez-Martín et al. (2014). □

Labbé et al. (2004) proposed using the *chain barrier constraints* (see Laporte et al. (1986)) or *path inequalities* (see Fischetti et al. (1998b)) for the Plant-cycle location problem, which has some similarity with this model from the network structure perspective.

**Proposition 5.** *For each $S \subseteq V, \{i, i'\} \subseteq I$ and a path $P := \{(i, i_1), (i_1, i_2), \ldots, (i_k, i')\}$ from $i$ to $i'$, the inequality:*

$$r(P) + \sum_{j \in S} z_{ij} + \sum_{j \notin S} z_{i'j} \leq |P| + 1 \qquad (24)$$

*is valid for BCCHRP polytope. Constraints* (8) *are special case of* (24) *when $P := \{[i, i']\}$.*

**Proposition 6.** *For all $S \subseteq V$*

$$r(\delta^+(S)) \geq \left\lceil \frac{\sum_{i \in S} \sum_{j \in V} z_{ij} + r(\delta^+(S))}{U} \right\rceil - \sum_{i \in S} z_{ii} \qquad (25)$$

$$r(\delta^-(S)) \geq \left\lceil \frac{\sum_{i \in S} \sum_{j \in V} z_{ij} + r(\delta^-(S))}{U} \right\rceil - \sum_{i \in S} z_{ii} \qquad (26)$$

*are valid inequalities.*

*Proof.* Given that $r(\delta(S)) = r(\delta^+(S)) + r(\delta^-(S)) = 2r(\delta^+(S)) = 2r(\delta^-(S))$, c.f. Proposition 3 in Rodríguez-Martín et al. (2014). □

**Proposition 7.** *For all $(i, i') \in A$ and $S \subseteq V$ where $i \in S$, $i' \in V/S$,*

$$r_{ij} + r_{ji} \leq \sum_{j \in V/S} z_{ij} + \sum_{j' \in S} z_{i'j'} \qquad (27)$$

*is a valid inequality.*

This inequality basically says that if for a given arc, its tail is in $S$ and its head is in $V/S$ (or vice versa), either the tail is allocated to a hub node in $V/S$ or the head is allocated to a hub node in $S$.

For every hub node $k$, the total flow entering to the hub (from the nodes not allocated to it) and travels to the first spoke on its route is the sum of demands of all the spoke nodes on the route minus demands of these spokes nodes that are fulfilled from the spoke nodes on the same route (the antecedent nodes supply part of the demands of the ensuing ones on the route). The latter is usually below a certain percentage of the whole demand of every node. A worst-case upper bound can be approximated given the number of nodes in the instance, the maximum possible number of nodes on a route and $max\{w_{ij} : \forall i, j\}$ in the flow matrix. We use $\beta$ to represent such approximations. This analogously holds for the case of the last arc visited on the route just before arriving to the hub node and the volume of flow on it.

**Proposition 8.** *For all $i \in V$, we define $O_i = \sum_{j \varsubsetneq i} w_{ij}$ and $D_i = \sum_{j \varsubsetneq i} w_{ji}$. Then, the following inequalities are valid for (BCCHRP):*

$$\sum_i O_i z_{ik} \leq \beta C z_{kk} \qquad (28)$$

$$\sum_i D_i z_{ik} \leq \beta C z_{kk} \qquad (29)$$

It must be noted that such valid inequalities are very often dominated by the capacity constraints (16) but as they are written in a way independent of the flow variables, they are mainly of interest in a Benders-like decomposition approach. Using such valid inequalities, we can filter out some of the solutions from the MP feasible space that render the Benders subproblem infeasible. This should conceivably contribute in a faster convergence of the algorithm.

**Proposition 9.** *Let $k$ be a hub node and $i$ be a spoke node allocated to it. The following inequalities are valid for (BCCHRP):*

$$r_{ki} D_i \leq C z_{kk} \qquad (30)$$
$$r_{ik} O_i \leq C z_{kk} \qquad (31)$$

*Proof.* If an arc $(k, i)$ is established between a hub node $k$ and a spoke node $i$ then $D_i$ cannot be larger than the capacity of the vehicle. Similarly, if the arc $(i, k)$ is established between a spoke node $i$ and a hub node $k$ then $O_i$ cannot be larger than the capacity of the vehicle. In other words, such nodes are candidates for being hub node rather than spoke node. □

Gendreau et al. (1998) proposed several classes of valid inequalities for the undirected Selective TSP (STSP). The following constraints can be easily derived from among those valid inequalities.

**Proposition 10.** *The following inequalities are valid for (BCCHRP):*

$$\sum_{j \in S : j \neq i} (r_{ij} + r_{ji}) \leq |S|, \qquad \forall i \in S \qquad (32)$$

Note that constraints (9) are special cases of constraints (32) where $|S| = 1$.

In Labbé et al. (2004) and Gendreau et al. (1998) also propose variants of the classical 2-matching constraints for a problem with some similarities in structure. These valid inequalities are also valid for our problem.

**Proposition 11.** *The following 2-matching constraints are valid for (BCCHRP):*

$$2 \sum_{i,j \in S} r_{ij} + \sum_{i \in S, j \in V/S} (r_{ij} + r_{ji}) \leq |S| + \frac{|A'| - 1}{2}, \qquad (33)$$

*for all $S \subset V$ and all $A' \subset A$ satisfying,*
*i)* $|\{i, j\} \cap S| = 1, \quad \forall (i,j) \in A'$,
*ii)* $|\{i, j\} \cap \{k, l\}| = \emptyset, (i, j) \neq (k, l) \in A'$, and
*iii)* $|A'| \geq 3$ and odd.

*Proof.* We have,

$$2 \sum_{i,j \in S} r_{ij} + \sum_{i \in S, j \in V/S} (r_{ij} + r_{ji}) = 2|S| \qquad (34)$$

and also since $r_{ij} \leq 1, \forall i, j \in E'$, by adding these inequalities to (34), taking into account that the left-hand side in (34) is an integer value, $|E'| \geq 3$ and $|E'|$ is odd, by rounding down the right-hand side, the proof is complete. □

## 3. Solution method

Benders decomposition (Benders, 2005, 1962) is a primal decomposition method, which has proven to be an efficient method in dealing with large-scale MIP models in facility location-type problems. The idea of Benders decomposition relies on relaxing the complicating variables (e.g. those variables that are integer) from the model. It then exploits the primal/dual relationship to generate cuts to separate the solutions of the master problem and tighten the outer approximation, until it proves optimality.

If a solution to the Benders master problem is not a feasible (partial) solution to the problem as a whole, the generated cut(s) in Benders approach is(are) feasibility cut(s) otherwise, we separate optimality Benders cut(s). It is well-known that feasibility cuts do not significantly contribute to the improvement of dual bound, and their role is to prevent infeasible solutions from the master problem. Given that, some authors, e.g. Gelareh and Nickel (2011) proposed replacing the natural choice of the master problem by an auxiliary model introduced earlier in Maculan et al. (2003) that guarantees connected subgraphs of a given set of nodes. Together with some branching rules, proposing a way of generating relative interior points and an appropriate decomposition of flow problem, this Benders decomposition based algorithm has shown to be very efficient and capable of solving realistic size instances in a very reasonable time. However, a trade-off between making the master problem harder to solve but tightening the MP polytope and eliminating the infeasible solution is not very trivial. On the other hand, Codato and Fischetti (2006) and Fischetti et al. (2009) proposed other approaches based on the *minimal infeasible subsystem* to deal with the infeasible solutions.

Several techniques have been proposed in the literature to improve performance of the Benders algorithm: Magnanti and Wong (1981a) deal with degenerate subproblem and choosing the best solution among many dual solutions; Magnanti et al. (1986) address the network design problems, in particular. Among many other contributions, we refer to Poojari and Beasley (2009), McDaniel and Devine (1977), and Sherali and Fraticelli (2002) for further reading. Benders decomposition for variants of hub location problems include de Camargo et al. (2009b, 2008) for models with emphasis on economies of scale and Contreras et al. (2011) for a classical model for the uncapacitated hub location problem.

The general framework of our algorithm is a Benders decomposition algorithm that is equipped with a cuts separation phase in addition to the usual Benders cuts. Such cuts are separated from among a set of valid inequalities introduced previously. The algorithm starts with a heuristic phase within which a relatively good feasible solution of the problem is found (if any), using an inexpensive local search algorithm. The Benders master problem is composed of binary variables (and the variables estimating the subproblem objective function) and the constraints uniquely containing those variables. Using such a feasible solution of the heuristic algorithm, the Benders subproblem can generate optimality cut(s) that will tighten –often to a great extent– the master problem in the very beginning of algorithm. Such cut(s) will most probably reduce the number of calls made to separation oracles (for valid inequalities as well as Benders cuts) dramatically. When separating Benders cuts, we use the techniques that make use of relative interior point of master problem allowing to choose the sharpest cuts. As long as a single non-Benders cut can be separated, we will not separate any Benders cut as such solutions –if representing an infeasible network structure– will only lead to generation of feasibility cuts that will cut the current recession ray and do not necessarily improve the lower bound significantly. Otherwise, we will separate Benders cuts. In generating Benders cuts, we use the method in Magnanti and Wong (1981a) to tighten our cuts. For this, we needed some interior points in the master polytope. The latter is approximated by the *pseudo-average* method in subsection 3.5. To do so, we ignore the capacity constraints and generate such points. Our numerical experiments show that any one of the methods proposed in subsection 3.5 produces cuts improving significantly the overall convergence when used in the Magnanti-Wong's technique. However, *pseudo-average interior* helps producing better cuts from the Magnanti-Wong's technique when compared to other two ways. In order to accelerate convergence, on a given frequency (or when certain conditions are fulfilled) we carry out an elimination test allowing to fix some variables and as a consequence getting rid of several dependent variables.

*3.1. Heuristic solution*

Even having an optimal solution at hand and adding the corresponding valid inequalities to the MP does not guarantee immediate convergence of this algorithm. In some special cases where the model is very tight, this might indeed occur. However, some authors –e.g. Poojari and Beasley (2009)– have proposed to use meta-heuristics to generate feasible solutions and add the corresponding cuts to the MP in order to accelerate convergence. Such an effort very often contributes to the bounds reduction on Benders subproblem cost estimation variable(s) $\eta$ (or $\eta_{ij}$) in the very early iterations and can sometimes reduce computational time significantly.

The general sketch of this heuristic is as follows: An initial solution is created by sorting the nodes with respect to the total (supply plus demand) flows. For every $l : q \leq l \leq p$, we do as follows. We choose the first $l$ nodes with the highest total volumes and designate them as hub nodes. Then using the same rule we allocate the remaining nodes to the hub nodes while checking the capacity constraints. For example, let $p = 3$ and $\{3, 4, 2, 1, 5, 7, 10, 6, 8, 9\}$ be the sorted list. We will have an initial solution as $S^0 = \{(3, 1, 10, 9), (4, 5, 6), (2, 7, 8)\}$. The elements of this set are themselves some ordered sets. E.g. $(2, 7, 8)$ represents a route with 2 representing the hub node and the spoke nodes 7 and 8 forming a route $2 \to 7 \to 8 \to 2$). Such a way of generating an initial solution guarantees a certain degree of fairness in distribution of volumes among different routes and it becomes more likely that the resulting solution respects capacity constraints —mainly when the fleet of vessels is homogeneous. We start from $S^0$ and for every node, in order, we try to find a better place by relocating it from a subset to another or perform a reordering within the same subset. The best move is selected, and the best-found solution ($S^0$) is updated. The next iteration starts from this new solution and applies the same relocation-reordering moves until a termination criterion is met. Our termination criteria is set to a count of $|V|$ non-improving iterations. Note that the computation of objective function is trivial once the network structure is given because there is a unique path from every origin to every destination.

It must be emphasized that this heuristic algorithm is not intended to be an independently sophisticated solution algorithm producing near optimal solutions. Rather, it suffices to have an *inexpensive* greedy procedure producing *relatively good* solutions (at least more meaningful than the initial solution reported by CPLEX for our master problem). The reason why we insist on such an inexpensive procedure is that as we do not have an off-line full description of our MP polytope (because of branch-and-cut scheme) even an available optimal solution might not guarantee an immediate convergence of Benders (perhaps not even in a few iterations) and investing too much in this procedure, in general, does not pay off.

It must be noted that, the master problem does not *directly* depend on any of the design variables and in the early iterations the MP resembles a feasibility problem because $\eta$ ($\eta_{ij}$) have not yet reached a good approximation and very limited information is available to help guiding towards a good solution and generating better cuts.

### 3.2. Initial relaxation

The initial relaxation of the model is a composed of a subset of constraints of BCCHRP with constraints (2)-(9). Some authors, for instance Gendreau et al. (1998) and Labbé et al. (2004) decided to also remove constraint (8) from the relaxation —perhaps due to the symmetry imposed by tournament constraints (9). However, in the case of this problem, we did not find any significant changes in numerical behavior in the presence or absence of such constraints.

The flow constraints (10)-(14) are handled using Benders decomposition. The capacity constraints (16) are also handled in the subproblem. Moreover, constraints (28)-(31) are added to the master problem in order to eliminate parts of the infeasibility. Constraints (28) and (29) act as a relaxed version of the original capacity constraints. The resulting Benders Master Problem (MP) follows:

$$min \ \eta = \sum_{ij} \eta_{ij} \tag{35}$$

$$s.\ t.$$
$$(2), (3), (4), (5),$$
$$(6), (7), (8), (9),$$
$$(28), (29), (30), (31),$$
$$r \in \mathbb{B}^{n^3}, z \in \mathbb{B}^{n \times n}, \eta_j \geq 0, \ \forall \ i, j \tag{36}$$

*Separation of Benders cuts.* There are some key aspects in generating Benders cuts for this problem. A Benders cut generated using (10)-(16) contains only $z_{ii}$ and $r_{ij}$ variables. All the allocation variables for the spoke nodes are absent in these cuts. On one hand, it is often better that every variable appears in one or more Benders cut(s) and total absence of a given subset of variables is not the best idea. On the other hand, we seek cuts that are not dense (the cuts become too similar from one iteration to another and causing numerical difficulties due to the larger eigen values of the LP matrix) and include a moderate size subset of variables in them.

In model, (10)-(16) we can replace some of the constraints as follows:

1. For constraints (10)-(11):

$$\sum_k z_{ik} = 1 \Rightarrow \sum_{k \subsetneq i} (x_{ijik} + s_{ijik}) = \sum_k z_{ik}, \qquad \forall i, j \subsetneq i, \tag{37}$$

$$\sum_k z_{ik} = 1 \Rightarrow \sum_{l \subsetneq j} (x_{ijlj} + s_{ijlj}) = \sum_k z_{ik}, \qquad \forall i, j \subsetneq i, \tag{38}$$

2. For constraints (13)-(14):

$$\sum_{k} \bar{z}_{ik} = 1 \Rightarrow \sum_{l \varsubsetneq k} \bar{x}_{ijkl} \leq 1 - \sum_{ik:k \varsubsetneq i} \bar{z}_{ik} \qquad \forall i,j,k : j \varsubsetneq i, k < l \qquad (39)$$

$$\sum_{k} \bar{z}_{ik} = 1 \Rightarrow \sum_{l \varsubsetneq k} \bar{x}_{ijlk} \leq 1 - \sum_{ik:k \varsubsetneq i} \bar{z}_{ik} \qquad \forall i,j,k : j \varsubsetneq i, k < l \qquad (40)$$

In this way, one can also accommodate variables $z_{ik}$, $\forall k \varsubsetneq i$ in the separated Benders cuts.

Let $\eta$ (or equivalently $\sum_{i,j \varsubsetneq i} \eta_{ij}$) be the continuous variable underestimating the objective value and the only variable in the objective function of the master problem, an optimality Benders cut resembles:

$$\sum_{ij \varsubsetneq i}(u^1_{ij} + u^2_{ij}) + \sum_{k}(u^4_{k} + u^5_{k})z_{kk} + \sum_{i,j,k,l:l \varsubsetneq k} u^6_{ijkl} r_{kl} + \sum_{k} u^7_{k} Cz_{kk} + \sum_{i,j,k,l:l \varsubsetneq k}(u^{10}_{ijkl} + u^{11}_{ijkl}) \leq \eta \qquad (41)$$

which is lacking variables $z_{kl}$, $\forall k, l \varsubsetneq k$. The feasibility cuts are obtained by replacing $\eta$ with 0.

The transformation (37) and (38) generate cuts of the following form:

$$\sum_{ij \varsubsetneq i}(u^1_{ij} + u^2_{ij})(z_{ij}) + \sum_{k}(u^4_{k} + u^5_{k})z_{kk} + \sum_{i,j,k,l:l \varsubsetneq k} u^6_{ijkl} r_{kl} + \sum_{k} u^7_{k} Cz_{kk} + \sum_{i,j,k,l:l \varsubsetneq k}(u^{10}_{ijkl} + u^{11}_{ijkl}) \leq \eta. \qquad (42)$$

Yet, given that $z_{kk} = 1 - \sum_{l \varsubsetneq k} z_{kl}$, the transformation **(39)** and **(40)** generate cuts of the form:

$$\sum_{ij \varsubsetneq i}(u^1_{ij} + u^2_{ij})(z_{ij}) - \sum_{k}(u^4_{k} + u^5_{k})(\sum_{l \varsubsetneq k} z_{kl}) + \sum_{i,j,k,l:l \varsubsetneq k} u^6_{ijkl} r_{kl} - \sum_{k} u^7_{k} C(\sum_{l \varsubsetneq k} z_{kl})$$
$$+ \sum_{i,j,k,l:l \varsubsetneq k}(u^{10}_{ijkl} + u^{11}_{ijkl}) + \sum_{k}(u^4_{k} + u^5_{k}) + \sum_{k} u^7_{k} C \leq \eta \qquad \mathbf{(43)}$$

that does not contain any $z_{kk}$ variable.

As a matter of fact, our numerical experiments suggest separating all three kinds of cuts (no one is traded for another). There is a small, yet important fact about the constraints (41)-(43). In (41), all the variable coefficients have the same sign (negative) but in (42) and (43) we have a mixture of pos/neg coefficients. From a computational perspective, such cuts are not contributing very well, at least when dealing with Benders decomposition.

An advantage of such transformed constraints is that one single Benders cut does not have to include all MP variables, and we can avoid dense constraints by simply distributing variables among different Benders cuts (41)-(43). Nevertheless, when all three forms of cut were added, the risk of numerical instabilities and some signs of singularity (in numerical sense) in the basis increases. Therefore, we choose to add the cuts by selecting their type according to a uniform distribution.

*Multi-cut Benders decomposition.* When a feasible solution to the original problem is available, the capacity constraints are either slack in which case the dual values are 0 (complementary slackness) or one can still move over the face of optimality and find alternative optimal solutions for which the dual values of capacity constraints are zero.

When $u^7_k = 0$, $\forall k$, the corresponding terms in the Benders cut vanish. As a result, one can split the Benders cut into $n(n-1)$ sub-cuts (corresponding to every $\eta_{ij}$). Such disaggregated cuts have shown to be very effective and can dramatically accelerate the lower bound improvement.

*Separation of valid inequalities* (18) *and* (19). The separation logic follows the fact that if a node $i$ is allocated to a hub node $j$ then if we send one unit of flow from a dummy node $s$ to the wherever $i$ is allocated, there must be finally one unit of flow that arrives to $i$. If such flow cannot be 1, then a cut is identified.

Let $s$ be a dummy node and establish a graph $G'(V \cup \{s\}, A')$ where there is an arc from $s$ to every $j$ if $\bar{z}_{ij} > 0$ with capacity $\bar{z}_{ij}$. Moreover, for every $\bar{r}_{ij} > 0$ we add an arc $(i,j)$ with capacity $\bar{r}_{ij}$. Now, a set $S \subset V$ where $i \in S, s \notin S$ defines a cut of $\delta^+(S)$ if the cut capacity is less than one. This set $S$ can be examined to identify both constraints (18) and (19). The complexity of this separation is based on the complexity of the max-flow algorithm of Boykov and Kolmogorov (2001). We opted to use Boykov-Kolmogorov max-flow algorithm as in graph $G'$, for every arc the reverse arc might be non-zero.

This is very similar to the procedure proposed in Rodríguez-Martín et al. (2014) and Labbé et al. (2004), except that in graph $G'$, the arc $(s,j)$ has capacity $\bar{z}_{ij}$ rather than $2\bar{z}_{ij}$.

*Separation of valid inequalities* (20) *and* (21). We establish a graph $G''(V, A'')$ where $A''$ is composed of arcs $(i,l)$ for which $r_{il} + z_{i'l} > 0$, arcs $(i',l)$ for which $r_{i'l} + z_{il} > 0$ and all other arcs $(k,l)$ where $r_{kl} > 0$.

A set $S \subset V$ where $i \in S, j \in V/S$ defines a cut if $\delta^+(S) < 2$. This set $S$ can be examined to identify both constraints (20) and (21). The complexity of this separation is based on the complexity of the max-flow algorithm of Boykov and Kolmogorov (2001). We opted to use Boykov-Kolmogorov max-flow algorithm as in $G'$ for some positive arcs the reverse arc might be non-zero.
Rodríguez-Martín et al. (2014) and Labbé et al. (2004) used a similar procedure for separating similar constraints.

*Separation of valid inequalities* (23). According to Rodríguez-Martín et al. (2014), such violated inequalities can be identified in a cut set $S$ of the support graph $G'''(V, A''')$ or in $V/S$. As in the separation procedure of constraints (18) (19), (20) and (21), we have already identified a cut-set $S$, one can also search for possible violation of (23) within these cut-sets or their complements. The complexity of this separation is based on the complexity of the max-flow algorithm of Boykov and Kolmogorov (2001).

*Separation of valid inequalities* (24). Separation of such inequalities is not an easy task. Moreover, our numerical experiments replicate similar observations as in Labbé et al. (2004). The effort invested in the separation of inequalities for $|P| \geq 2$ does not pay off. Therefore, we stick to the constraints (8) that already exist in the initial relaxation.

*Separation of valid inequalities* (25) *and* (26). Let $U$ be the maximum possible arcs along a given route. We establish a graph $G''''(V \cup \{s,t\}, A'''')$ where $A''''$ is composed of arcs $(i,j)$ for which $(U-1)\bar{r}_{ij} > 0$ with a capacity $(U-1)\bar{r}_{ij}$, all arcs $(s,i)$ for which $U\bar{z}_{ii} > 0$ with capacity $U\bar{z}_{ii}$ and all arcs $(i,t)$ with a capacity $\sum_{j \in V} \bar{z}_{ij}$. A set $S \subset V$ where $t \in S, s \in V/S$ defines a violated cut, if the cut capacity is smaller than $|V|$ (see Rodríguez-Martín et al. (2014) for a similar separation).

*Separation of valid inequalities* (27). The separation of such constraints is exactly the same as the one in Rodríguez-Martín et al. (2014). Let $S = \{i\} \cup \{j \in V - \{i'\} : z_{ij} \geq z_{i'j}\}$, for a given arc $(i,i')$. According to Rodríguez-Martín et al. (2014), if the inequality is not violated for this choice of $S$, it is neither violated for any other set $S'$. The complexity is $O(|V|^3)$.

*Separation of valid inequalities* (32). The separation of this class of valid inequalities is reduced to solving a Knapsack Problem (KP) with two additional constraints. The weights of items correspond to $\bar{r}_{ij}$.

$$\max \sum_i \bar{b}_i$$
$$s.t. \sum_{i,j} \bar{r}_{ij} s_{ij} \geq \sum_i \bar{b}_i + Q$$

$$s_{ij} \leq b_i \qquad \forall i, j$$
$$s_{ij} \leq b_j \qquad \forall i, j$$
$$b_i, s_{i,i} \in \{0, 1\} \tag{44}$$

Any feasible solution to this problem identifies a set $S = \{i : b_i = 1, \forall i\}$ that is used to separate the corresponding cut. An alternative separation algorithm was proposed in Gendreau et al. (1998).

*Separation of valid inequalities* (33). Fischetti et al. (1998a) proposed a heuristic separation algorithm for such constraints. It suffices to establish the support graph and use a threshold value to remove the edges for which the values of $r_{ij}$ are more than the threshold and seek for $S$ among the connected components of the reduced graph.

### 3.3. Reduction test

In our problem, none of the binary design variables of master problem has a cost and the whole cost is due to the routing part, which remains in the subproblem. The reduced costs associated to every non-basic variable is directly related to the coefficients of the generated Benders cuts. However, this is quite tricky as such coefficients of variables in a Benders cut represent overestimates of reduced costs of those variables. Therefore, they are not sufficiently reliable to be exploited for any kind of elimination test.

In the Benders subproblem we have multi-commodity flow problem with capacities on the non-hub arcs. Let $\bar{r}_{ij}, \bar{z}_{ij}$ be the solution to the Benders master problem with a feasible network structure (even in the fractional sense).

$$\min \sum_{i,j,k,l} t_{kl}(s_{ijkl} + \alpha x_{ijkl}) + \sum_{ijkl: \text{all distinct}} (\phi^k + \phi^l) x_{ijkl} \tag{45}$$

s. t.

$$\mathbf{u}^1 : \sum_{k \neq i} (x_{ijik} + s_{ijik}) = 1, \qquad \forall i, j \in V : j \neq i, \tag{46}$$

$$\mathbf{u}^2 : \sum_{l \neq j} (x_{ijlj} + s_{ijlj}) = 1, \qquad \forall i, j \in V : j \neq i, \tag{47}$$

$$\mathbf{u}^3 : \sum_{l \neq i,k} (x_{ijkl} + s_{ijkl}) = \sum_{l \neq j,k} (x_{ijlk} + s_{ijlk}), \qquad \forall i, j, k \in V, k \notin \{i, j,\}, \tag{48}$$

$$\mathbf{u}^4 : \sum_{l \neq k} x_{ijkl} \leq z_{kk} \qquad \forall i, j, k \in V : j \neq i, k < l \tag{49}$$

$$\mathbf{u}^5 : \sum_{l \neq k} x_{ijlk} \leq z_{kk} \qquad \forall i, j, k \in V : j \neq i, k < l \tag{50}$$

$$\mathbf{u}^6 : s_{ijkl} \leq r_{kl} \qquad \forall i, j, k, l \in V : l \neq k \tag{51}$$

$$\mathbf{u}^7 : \sum_{ijl: j \neq i} w_{ij} s_{ijkl} \leq C, \qquad \forall k \in V \tag{52}$$

$$\mathbf{u}^8 : r_{ij} = \bar{r}_{ij}, \qquad \forall k \in V \tag{53}$$

$$\mathbf{u}^9 : z_{ij} = \bar{z}_{ij}, \qquad \forall k \in V \tag{54}$$

$$\mathbf{u}^{10} : s_{ijkl} \leq 1, \qquad \forall k \in V \tag{55}$$

$$\mathbf{u}^{11} : x_{ijkl} \leq 1, \qquad \forall k \in V \tag{56}$$

$$r \in \mathbb{B}^{|V|^3}, z \in \mathbb{B}^{|V| \times |V|}, x_{ijkl}, s_{ijkl} \in \mathbb{R}_{[0,1]}^{|V|^4} \tag{57}$$

Our numerical experiments reveal that, as expected, constraints (53)-(54) may result in numerical instability. The remedy is to convert them to half-space constraints (i.e. $\leq$) without any harm as the objective is a minimization problem. For every set of constraints in this problem the corresponding dual vector is indicated on the left side of the constraints (see $\mathbf{u}^i$). $\mathbf{u}^8$ and $\mathbf{u}^9$ correspond to the correct reduced costs of the design variables. These costs can be used both in generating Benders cuts and in variable fixing (reduction test).

Given a lower bound and an upper bound on the objective value, local and global cuts can be added to the branch-and-bound nodes. At the root node, for a given variable $v$, if $|rc(v)| > UB - LB$ then the variable is set to its current value and $v = \bar{v}$ is a globally valid cut. On any other node, such a fixing is a local cut.

Furthermore, we can do the following:

a) Let $\eta$, that is the linking variable of the master problem underestimating the real objective function, be re-written as $\sum_{i,j \in C_i} \eta_{ij}$, $\eta_{ij} \geq 0$, $\forall i, j$ where $\eta_{ij}$ underestimates the $(i, j)$-flow. One can reduce bounds of $\eta_{ij}$ variables. From *below*, $\eta_{ij}$ cannot be less than $\alpha t_{ij}$ because the minimum is taken when $i$ and $j$ are directly connected and both are hub nodes. From *above*, we can have at most $n - (p - 2) \times q - 2$ arcs (not counting the hub-to-hub edge) along any O-D pair. Therefore, if one wishes to avoid solving a *constrained longest path problem*, one needs to find the $n - (p - 2) \times q - 2$ most expensive arcs, consider two most expensive transshipment and set it as an upper bound on $\eta_{ij}$. However, there is not much to exploit in terms of numerical efficiency using such an upper bound.

b) As soon as a feasible solution to the problem is known, both lower (LB) and upper bounds (UB) are known. We also have dual values associated with the subproblem constraints of Benders algorithm and a list of non-basic design variables, which are subject to our reduction test. We do the following:
  - *cover cuts or combinatorial cuts*: For any minimal number of variables in the set $X$ for which the sum of their reduced costs is larger than $LB - UB$ we add a cover cut $\sum_{i=1}^{|X|} x_i \leq |X| - 1$. This sort of reduction is useful in the very early iterations of Benders as the gap is large and it is rather hard to exploit primal/dual information and single out a variable to be fixed. In some cases, small IP with a few constraints needs to be solved to find a meaningful set (such constraints include those for avoiding arcs in the opposite direction (tournaments)).
  - *single fixing*: Later in the iterations, as the optimality gap reduces, finding a single variable to be fixed becomes easier. If the corresponding reduced cost does not suggest a variable to take 1 in any optimal solution, we can safely fix it to 0.

## 3.4. Branch-and-cut algorithm

As mentioned earlier, our emphasis in our branch-and-bound is on fast lower bound improvements. That is, the nodes are explored using breadth-first scheme. For that, we separate all possible violated cuts at the root node. If no inequality from (18)-(27) is separated then we generate Benders cuts.

The Benders subproblem is comprised of constraints (10)-(16). When solving Benders subproblem, we generate the cuts based on the idea of non-dominated cuts proposed by Magnanti and Wong (1981b) while instead of the computing a relative interior point, we resort to the methods proposed in subsection 3.5. When separating Benders cuts, we do not separate any cut for the fractional solution, if the relative gap is less than 10%.

We have chosen to give priority to improving dual bound in our branch-and-bound procedure, and also decided to branch on $z$ variables first. It might look more evident to branch on $x_{ij}$ because they have a direct impact on the status of $z_{ik}$. However, the numerical experiments do not confirm this. Moreover, to have a more effective variable elimination/fixing procedure, we need to quickly obtain a reasonable evaluation of $z_{ik}$ reduced costs such that our reduction test performs more effectively and variables are fixed/removed in the early iterations, if any.

In addition to the above rules, occasionally, we test another possibility for branching. Let $G'(V, A')$ be the support graph of any given feasible solution to the problem that is obtained by removing all the hub-level interconnections.

**Lemma 1.** *Any support graph $G'$ resulted by a integer feasible solution is composed of at least $q$ and at most $p$ strongly connected.*

Note that based on the preceding lemma, the cut size of every integer feasible solution is 0.

In the solution process, we verify at a given frequency whether the support graph contains a correct number of components. If not, then we run $|V|(|V|-1)/2$ min-cut algorithms on the directed graph $G'$. If in a feasible (fractional) solution we end up having a cut size of at least $1-\varrho, |\varrho| \leq 0.05$, we have our candidate (a variable or a set of variables) to branch on.

*3.5. Enhancement to the classical Benders decomposition*

Some key tunings are needed in order to better exploit the efficiency of Benders decomposition, particularly when the problem embodies a multicommodity flow problem (due to primal degeneracy/dual multiple optimality), which is the case in variants of HLPs. Interested readers are referred to Contreras et al. (2011) and de Camargo et al. (2009b, 2008) among others.

Magnanti and Wong (1981a) proposed using relative interior point to generate non-dominated cuts. This technique has shown to be very effective and has reported very promising results in the literature. This includes Gelareh and Nickel (2011), wherein the authors experienced a dramatic improvement in convergence. Nevertheless, computing such a relative interior point is not always an easy task. Although Freund et al. (1985) propose an alternative LP model for computing a relative interior point, however, the proposed model is equivalent to solving an LP relaxation of the master problem each time a new relative interior is needed. Taking into account that the size of the master problem is getting larger each time a cut (no matter optimality or feasibility) is added, computing such a node becomes more and more expensive. This is particularly the case when more than one cut (e.g. disaggregated Benders cuts) is added to the master problem at each iteration.

**Definition 1.** *Let represent a solution $C^H$ to this problem as $C^H = (C^1, C^2, \ldots, C^h)$ where $C^i$ is the route associated to the hub $i$. Such a $C^h$ is characterized by a sequence of spoke nodes traversed on the route. Let $\overleftarrow{C^h} := \overrightarrow{C^h}$ and define $\overleftarrow{C^i}$ to be the same route $C^i$ where the arcs are in the opposite directions (i.e. reverse route). Let also $f(g)$ be a function where $f^0(\overrightarrow{C^i}) = \overrightarrow{C^i}$ and $f^1(\overrightarrow{C^i}) = \overleftarrow{C^i}$. A mapping $F : [\Pi]_{1 \times h} \mapsto (f^{\Pi_1}(C^1), f^{\Pi_2}(C^2), \ldots, f^{\Pi_h}(C^h))$ where $\Pi_i \in \{0, 1\}, \forall i$ is a bijective map that in the absence of capacity constraints, produces symmetric solutions having the same objective values. Each one of such solutions is differentiated from $C^H$ only with respect to the direction of routes associated to those non-zero entries of $[\Pi]_{1 \times p}$ (see Figure 3). Note that, $F : (0, 0, \ldots, 0) \mapsto (\overrightarrow{C^1}, \overrightarrow{C^2}, \ldots, \overrightarrow{C^h})$ and $F : (1, 1, \ldots, 1) \mapsto (\overleftarrow{C^1}, \overleftarrow{C^2}, \ldots, \overleftarrow{C^h})$.*

**Lemma 2.** *For any solution $C^H$, in the absence of capacity constraints, $\overline{C^H} = \frac{\sum_{i=2} w_i F(\Pi_i(C^H))}{\sum_i w_i}$, for distinct vectors $\Pi$,*

- *is a feasible solution,*
- *is a relative interior point of the polyhedron.*
- *if $C^H$ is an optimal solution, $\overline{C^H}$ can be chosen to be a non-corner point on the same optimality face.*

In the absence of the capacity constraints (16), three techniques are proposed to generate additional points (see Figure 2):

**pseudo-average interior** Here we rely on Lemma 2 with $w_i = 1 : i \in S, |S| > 1$. This gives an interior point (see $\hat{\mathbf{y}}$ in Figure 2). If all such symmetric solutions could be used in Lemma 2, such an interior point is indeed the real center of the face (see Figure 2 where $\mathbf{x}$ is the optimal, $p_7, p_8, p_{16}, p_{20}$ and $p_{19}$ are on the same optimality face and $\hat{\mathbf{y}}$ is the center of $p_7, p_8$ and $p_{20}$. $\hat{\mathbf{x}}$ is the center of all the points on the face).

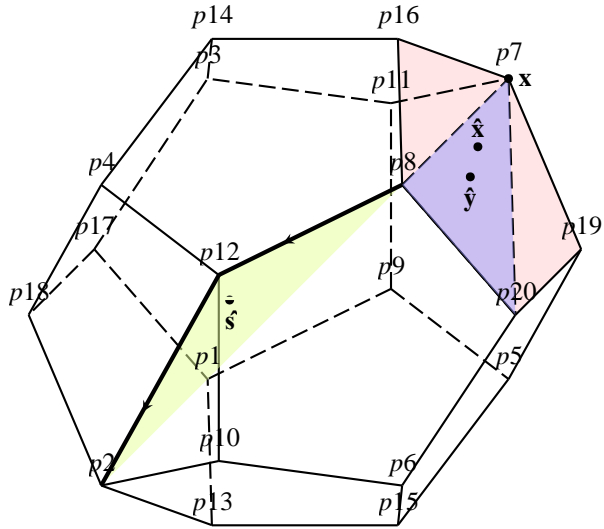

**Figure 2:** A hypothetical polytope and two ways of generating the non-corner interior points. Let this polytope represent the convex hull of integer feasible solutions: 1) Let $p8$, $p16$, $p7$, $p19$ and $p20$ be all on the same optimality face. I.e. they have the same objective but they are different in the direction of at least one route. The center of this optimality face is $\hat{\mathbf{x}}$. However, as enumerating such set of nodes in a higher dimension becomes intractable, one can restrict oneself to a subset of such nodes, namely $p8$, $p7$ and $p20$ in this case, and calculate a pseudo-center $\hat{\mathbf{y}}$. However, such a point can be chose to be non-corner point. 2) Let $p8$, $p12$ and $p2$ be the three most recent incumbents reported by the solver. The center of these corner points represented by $\hat{\mathbf{s}}$ is a node that *might be* inside the polytope, otherwise we try to converted to an interior point by moving along the line connecting this point to the current incumbent.

**uniform interior** Let $l = \frac{p+q}{2}$. We generate a fractional feasible solution as $z_{ii} = \frac{1}{n}$, $z_{ij} = \frac{1-n}{n-1}$ and $r_{ij} = \frac{n}{n(n-1)}$. Such a solution is an interior point of the polytope and is very close to 0 with respect to the Euclidian metric. Note that such a point is calculated once and never changed again. It is possible that capacity constraints or a Benders cut separates this point from the polytope in the very early iterations. However, when capacity is not very tight often this point remains within the polytope until after quite a few Benders iterations (as is the case here).

**incumbent relinking interior** In a branch-and-cut implementation of Benders decomposition, we store all the accepted incumbents in the course of the branch-and-bound and calculate their average (center). There are three possible situations: 1) Such a center is an interior point, 2) an interior point can be created from this by moving towards the interior of polytope in the direction of the line connecting this point to the current incumbent, 3) none of the first two cases holds and we fail to produce an interior solution. In Figure 2, let $p_2, p_{12}$ be the last two visited incumbents that may not currently part of the current polytope anymore and $p_8$ be the current incumbent. $\hat{\mathbf{s}}$ is the center of these point which is not guaranteed to be part of the polytope but it might be possible to obtain an interior point from it.

Our preliminary computational experiments suggest *pseudo-average* instead of others.

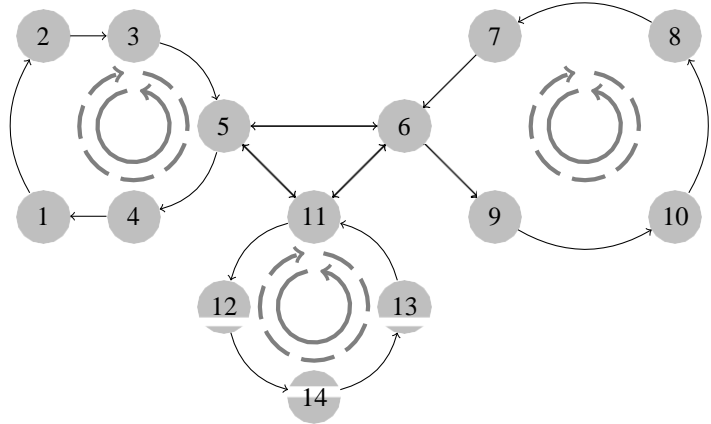

**Figure 3:** Symmetry in this problem.

## 4. Computational experiments

In this section, we elaborate on generating appropriate instances for our computational experiments. We make a comparison between the LP bounds of both models and present our computational results of applying the proposed decomposition on the instances of the testbed.

We have generated our instances based on the well-known Australian Post (AP) dataset (see Ernst and Krishnamoorthy (1999)). The transshipment times are randomly generated real values within [2, 5] for the given time unit. The original capacities from the AP instances are not used due to the feasibility issues. We assume that our fleet of vehicles is homogenous and therefore, the capacity is unique. Let $O_i = \sum_j w_{ij}$ and $D_j = \sum_i w_{ij}$. Let also $q$ and $p$ be the minimum and maximum number of hub nodes. The capacities are generated as in the following: $C_k^1 = \frac{(\sum_i O_i + \sum_i D_i)}{2 \times p}$. In our instances, $q$ is set to 3 in all the cases and $p \in \{3, 4, \ldots, \lfloor N/3 \rfloor\}$. The minimum number nodes allocated to a hub node (including itself), namely $\Gamma$, is set to 3.

### 4.1. Computational experiments with CPLEX

First, we used CPLEX 12.6.1 as the general-purpose MIP solver. A time limit of 36000 (10 hours) seconds has been set, and the number of threads is limited to 6. In Table 3, we report the results of our numerical experiments with the model (BCCHRP).

The first column reports the instant names in the format $nN\_pP\_\alpha$ where $N$ is the number of nodes, $P$ is the maximum number of hub nodes, and $\alpha$ represents the factor of economies of scale varying in $\{0.7, 0.8, 0.9\}$. The second column indicates CPU times elapsed during the LP resolution. The next column reports the LP objective function values followed by a column for the LP statuses. The next column report the value of $\sum_i z_{ii}$ in the LP solution. In columns 'TimeIP' and 'IPobj', we report the CPU times for IP resolution and the objective function values, respectively. The subsequent column reports the number of nodes processed before termination criteria is met. The column 'IPStatus' reports the termination criterion that has been met for every instance. The next column reports the optimality gaps when CPLEX decided to terminate.

In Table 3, one observes that the LP bound is rather weak. In the LP relaxation, there is a tendency towards opening maximum possible number of (partial) hubs (in the fractional sense). While for $n = 10$ the instances were solved in less than one hour, for $n = 15$, optimality could not be proven within 10 hours. Slight deviation from 10 hours in CPLEX is due to the fact that CPLEX terminates once the last node is solved and does not interrupt LP solution. As $p$ increases, the objective function value tends to decrease. There are some signs of potential numerical issues for $n = 10$ as CPLEX terminates upon meeting its internal tolerance-related termination criteria.

**Table 3:** Numerical experiments with the 2-index model.

| instance | TimeRoot | LPobj | LPStatus | $\sum_i z_{ii}$ | TimeIP | IPobj | Nnodes | IPStatus | Gap (%) |
|---|---|---|---|---|---|---|---|---|---|
| n10_p3_0.9 | 3 | 2219.32 | Optimal | 3 | 1973 | 3395.63 | 15905 | Optimal | 0.00 |
| n10_p3_0.8 | 2 | 2160.99 | Optimal | 3 | 2176 | 3315.81 | 19729 | OptimalTol | 0.01 |
| n10_p3_0.7 | 2 | 2097.88 | Optimal | 3 | 3371 | 3235.99 | 50926 | OptimalTol | 0.01 |
| n15_p3_0.9 | 8 | 5998.29 | Optimal | 3 | 36008 | 14461.60 | 2923 | AbortTimeLim | 56.71 |
| n15_p3_0.8 | 7 | 5863.59 | Optimal | 3 | 36016 | 14702.10 | 4661 | AbortTimeLim | 56.97 |
| n15_p3_0.7 | 9 | 5699.74 | Optimal | 3 | 36016 | 11842.20 | 4305 | AbortTimeLim | 48.50 |
| n15_p4_0.9 | 24 | 5568.42 | Optimal | 4 | 36034 | 10785.20 | 5615 | AbortTimeLim | 43.16 |
| n15_p4_0.8 | 21 | 5419.25 | Optimal | 4 | 36029 | 9594.57 | 5072 | AbortTimeLim | 37.67 |
| n15_p4_0.7 | 15 | 5240.83 | Optimal | 4 | 36025 | 10425.50 | 6505 | AbortTimeLim | 44.79 |
| n15_p5_0.9 | 20 | 5203.81 | Optimal | 5 | 36020 | 9122.92 | 4391 | AbortTimeLim | 37.44 |
| n15_p5_0.8 | 16 | 5039.83 | Optimal | 5 | 36023 | 7381.18 | 7004 | AbortTimeLim | 23.58 |
| n15_p5_0.7 | 19 | 4852.73 | Optimal | 5 | 36026 | 7885.72 | 7149 | AbortTimeLim | 32.88 |

*4.2. Computational experiments with branch-and-cut algorithm*

The Branch-and-Cut algorithm makes use of CPLEX as a modern MIP solver offering several callbacks and functionalities to allow users to take over control of some parts of the whole solution process by implementing some customizations (interfering in branching, cutting, node resolution, etc.). A single thread has been used here and CPLEX emphasis is set to improve lower bound. When using CPLEX, one can implement both LazyConstraint and CutCallbacks to implement this algorithm. All the added cuts are added globally, and we avoid adding local cuts.

For small instances, the reduction test was very effective in eliminating up to around 60% percent of the variables at the root node. However, starting from $n = 15$, the CPU time elapsed in solving the auxiliary LP to obtain the dual values does not pay off and even such dual values were less useful as the LP relaxations become too weak. Therefore, we decided to not report any statistics in this regard in Table 4.

We have used the well-known Boost Graph Library implementations in our separation algorithms. For the second phase of our algorithm, the stopping criteria are, 1) a time limit of 7,200 seconds, 2) a gap of 5%, and 3) existence of a feasible solution incumbent. The termination criteria for the third phase are, 1) a time limit of 36,000 seconds (10 hours), 2) a gap of 5%, and 3) existence of a feasible solution incumbent. Table 4 reports the computational experiments with the Branch-and-Cut algorithm applied to the 2-index formulation.

The first column in the table indicates the instance name. The second column report CPU times in seconds. The next column reports the best incumbent found. The column '#Nodes' reports the number of nodes processed in the course of solution process. The column 'CplexStatus' reports the CPLEX status upon termination. The column 'Gap(%)' reports the termination gaps (the LB is also reported when the method failed). The next two columns represent number of feasibility Benders cuts ('#F. Cuts') and number of optimality Benders cuts ('#O. Cuts'). The next three columns report number of cuts separated from every class, namely constraints (20),(21) , (23) , (25) and (26). Some classes of valid inequalities were never generated and therefore, we removed the corresponding column. The columns 'H. O.' and 'H. T.' report the heuristic objective function values and the CPU times in seconds, respectively.

In general, the number of nodes processed is very moderate and most of the instances were solved to optimality. For two instances we had numerical issues that led to failures. The efforts to avoid such numerical instabilities by tuning different parameters and tolerances were not successful. In order to avoid numerical issues we have set the gap tolerance to 0.5%. One observes that the optimality is subject to the user-defined tolerance gap. As explained earlier, in the case of separating feasibility cuts, only one cut is separated and when optimality, multiple cuts can be added. Therefore, the number of feasibility Benders cuts corresponds to the iterations where infeasible solution has been encountered. These numbers are also very small. It means that the number of iterations with no (or very minimal) improvement in the lower bound is quite small. Regarding the separated non-Benders cuts, constraints (20),(21) and (23) are the most violated ones while very few constraints (25) and (26) are separated for some instances —when $n$ = 20 and $p$ = 6.

**Table 4:** Computational experiments with the Branch-and-Cut algorithm based on Benders decomposition.

| instance | TotalTime | Obj. Val. | #Nodes | CplexStatus | Gap(%) | #F. Cuts. | #O. Cuts. | (20),(21) | (23) | (25), (26) | H. O. | H. T. |
|---|---|---|---|---|---|---|---|---|---|---|---|---|
| n10_p3_0.9 | 14 | 3395 | 0 | Optimal | 0.37 | 0 | 696 | 2 | 2 | 0 | 3455.18 | 1 |
| n10_p3_0.8 | 14 | 3315 | 0 | Optimal | 0.42 | 0 | 617 | 2 | 0 | 0 | 3385.26 | 1 |
| n10_p3_0.7 | 17 | 3235 | 0 | OptimalTol | 0.09 | 0 | 785 | 0 | 0 | 0 | 3297.01 | 2 |
| n15_p3_0.9 | 398 | 11244 | 9 | OptimalTol | 0.19 | 3 | 3009 | 16 | 6 | 0 | 11391.28 | 4 |
| n15_p3_0.8 | 1395 | 11255 | 25 | OptimalTol | 0.27 | 3 | 6465 | 48 | 12 | 0 | 10840.70 | 4 |
| n15_p3_0.7 | 929 | 11120 | 47 | OptimalTol | 0.41 | 4 | 5681 | 14 | 20 | 0 | 11120.20 | 4 |
| n15_p4_0.9 | 899 | 9067 | 84 | OptimalTol | 0.00 | 1 | 2690 | 54 | 9 | 0 | 9067.25 | 9 |
| n15_p4_0.8 | 883 | 10170 | 39 | OptimalTol | 0.41 | 0 | 4706 | 14 | 2 | 0 | 10205.32 | 9 |
| n15_p4_0.7 | 1073 | 10683 | 97 | OptimalTol | 0.37 | 9 | 6039 | 18 | 0 | 0 | 10683.44 | 9 |
| n15_p5_0.9 | 1615 | 7387 | 231 | Optimal | 0.00 | 0 | 2364 | 54 | 61 | 0 | 7387.35 | 15 |
| n15_p5_0.8 | 591 | 7143 | 46 | Optimal | 0.00 | 0 | 1989 | 8 | 16 | 0 | 7143.31 | 13 |
| n15_p5_0.7 | 626 | 6899 | 68 | Optimal | 0.00 | 0 | 1581 | 12 | 16 | 0 | 6899.27 | 14 |
| n20_p3_0.9 | 28673 | 24935 | 288 | OptimalTol | 0.25 | 10 | 13250 | 290 | 119 | 0 | 25248.85 | 11 |
| n20_p3_0.8 | — | — | 114 | failed | 22.07%(LB = 19509.37) | — | — | — | — | — | 25035.00 | 12 |
| n20_p3_0.7 | 15307 | 5548 | 965 | OptimalTol | 0.46 | 33 | 747 | 526 | 206 | 0 | 24863.99 | 17 |
| n20_p4_0.9 | 15767 | 21585 | 184 | Optimal | 0.00 | 0 | 9202 | 188 | 77 | 0 | 21585.40 | 22 |
| n20_p4_0.8 | 14432 | 21284 | 107 | OptimalTol | 0.48 | 3 | 13217 | 88 | 35 | 0 | 21284.45 | 23 |
| n20_p4_0.7 | — | — | 199 | failed | 30.57%(LB = 15169.85) | — | — | — | — | — | 21850.60 | 23 |
| n20_p5_0.9 | 19051 | 17441 | 301 | Optimal | 0.00 | 0 | 6676 | 248 | 87 | 0 | 17441.27 | 53 |
| n20_p5_0.8 | 19572 | 17624 | 247 | OptimalTol | 0.24 | 16 | 5694 | 224 | 95 | 0 | 17624.71 | 33 |
| n20_p5_0.7 | 6470 | 16933 | 51 | Optimal | 0.00 | 0 | 5984 | 50 | 63 | 0 | 16933.63 | 37 |
| n20_p6_0.9 | 98991 | 17274 | 1525 | AbortUser | 31.33 | 0 | 10308 | 376 | 508 | 33 | 17274.12 | 44 |
| n20_p6_0.8 | 46695 | 16922 | 719 | AbortUser | 31.01 | 0 | 10462 | 138 | 77 | 1 | 16922.59 | 48 |
| n20_p6_0.7 | 27531 | 15738 | 275 | OptimalTol | 0.41 | 0 | 16339 | 52 | 50 | 0 | 16861.48 | 49 |

Our heuristic produced and improves the initial solution in a moderate time when compared to the total time spent in the whole algorithm.

## 5. Conclusions

In this paper, we have proposed a mathematical model for the Hub Location Routing Problem where to every hub node at least two spoke nodes are allocated forming a route starting from and terminating to the hub node itself. The routes must respect the capacity of the vehicle and the fleet of vehicles is assumed to be homogenous. The hub-level network is assumed to be fully interconnected. The problem is of interest in multi-modal distribution systems where there is a point-to-point connection between the major hubs while the local distributions takes place based on a given itinerary using a homogenous fleet of vehicles. We have identified several classes of valid inequalities (some of which were previously proposed for relevant problems in the literature) and have adapted the separation procedures. It has been shown that the state-of-the-art general-purpose MIP solvers such as CPLEX are not capable of dealing with even small size instances of the problem in reasonable time. Therefore decomposition-based methods need to be considered. A branch-and-cut algorithm based on a Benders decomposition has been proposed. In order to separate tighter Benders cuts, we have shown how symmetry can be exploited to extract interior solutions that can be used to generate non-dominated Benders cuts using Magnanti-Wong's method. We have suggested several enhancements to the classical Benders decomposition such as possibility of reducing the bounds of variables and fixing variables based on the reduced costs. Our computational experiments show that although our decomposition scheme is limited in the size of instances that can be solved, the optimal solutions are obtained with much less computational effort and at a fraction of time that a state-of-the-art MIP solvers would require.

Obviously, the success of many exact solution methods –particularly, the cutting-plane based ones– rely heavily on the polyhedral knowledge of the problem. The better the polytope of the problem is known, the better models can be developed and sharper cuts may be identified. Therefore, polyhedral study of this problem, heuristics and hybrid of exact and heuristic methods deserve more attention.

# References

Alumur, S., Kara, B. Y., 2008. Network hub location problems: the state of the art. European Journal of Operational Research 190, 1–21.

Alumur, S., Kara, B. Y., 2009. A hub covering network design problem for cargo applications in Turkey. Journal of the Operational Research Society 60, 1349–1359.

Alumur, S., Kara, B. Y., Karasan, O. E., 2009. The design of single allocation incomplete hub networks. Transportation Research Part B 43, 936–951.

Benders, J., 2005. Partitioning procedures for solving mixed-variables programming problems. Computational Management Science 2 (1), 3–19.

Benders, J. F., 1962. Partitioning procedures for solving mixed-variables programming problems. Numerische mathematik 4 (1), 238–252.

Boykov, Y., Kolmogorov, V., 2001. An experimental comparison of min-cut/max-flow algorithms for energy minimization in vision. In: Energy minimization methods in computer vision and pattern recognition. Springer, pp. 359–374.

Campbell, J. F., Ernst, A. T., Krishnamoorthy, M., 2002. Hub location problems. In: Drezner, Z., Hamacher, H. W. (Eds.), Facility Location: Applications and Theory. Springer, pp. 373–407.

Campbell, J. F., O'Kelly, M. E., 2012. Twenty-five years of hub location research. Transportation Science 46 (2), 153–169.

Cetiner, S., Sepil, C., Sural, H., 2006. Hubbing and routing in postal delivery systems. Tech. rep., Industrial Engineering Department, Middle East Technical University, 06532 Ankara, Turkey.

Çetiner, S., Sepil, C., Süral, H., 2010. Hubbing and routing in postal delivery systems. Annals of Operations Research 181 (1), 109–124.

Codato, G., Fischetti, M., 2006. Combinatorial Benders' cuts for mixed-integer linear programming. OPERATIONS RESEARCH-BALTIMORE THEN LINTHICUM- 54 (4), 756.

Contreras, I., Cordeau, J.-F., Laporte, G., 2011. Benders decomposition for large-scale uncapacitated hub location. Operations research 59 (6), 1477–1490.

Contreras, I., Fernández, E., Marín, A., 2009. Tight bounds from a path based formulation for the tree of hub location problem. Computers & Operations Research 36, 3117–3127.

Contreras, I., Fernández, E., Marín, A., 2010. The tree of hubs location problem. European Journal of Operational Research 202, 390–400.

de Camargo, R. S., de Miranda, G., Løkketangen, A., 2013. A new formulation and an exact approach for the many-to-many hub location-routing problem. Applied Mathematical Modelling 37 (12), 7465–7480.

de Camargo, R. S., de Miranda Jr., G., Ferreira, R. P. M., Luna, H. P. L., 2009a. Multiple allocation hub-and-spoke network design under hub congestion. Computers & OR 36 (12), 3097–3106.

de Camargo, R. S., de Miranda Jr., G., Luna, H. P. L., 2009b. Benders decomposition for hub location problems with economies of scale. Transportation Science 43 (1), 86–97.

de Camargo, R. S., Miranda, G., Luna, H., 2008. Benders decomposition for the uncapacitated multiple allocation hub location problem. Computers & Operations Research 35 (4), 1047–1064.

Ernst, A., Krishnamoorthy, M., 1999. Solution algorithms for the capacitated single allocation hub location problem. Annals of Operations Research 86 (0), 141–159.

Farahani, R. Z., Hekmatfar, M., Arabani, A. B., Nikbakhsh, E., Apr. 2013. Survey: Hub location problems: A review of models, classification, solution techniques, and applications. Comput. Ind. Eng. 64 (4), 1096–1109.

Fischetti, M., Gonzalez, J. J. S., Toth, P., Feb. 1998a. Solving the orienteering problem through branch-and-cut. INFORMS J. on Computing 10 (2), 133–148.

Fischetti, M., Gonzï¿½lez, J. J. S., Toth, P., 1998b. Solving the orienteering problem through branch-and-cut. INFORMS Journal on Computing 10, 133–148.

Fischetti, M., Salvagnin, D., Zanette, A., 2009. Minimal infeasible subsystems and Benders cuts. Mathematical Programming to appear.

Freund, R. M., Roundy, R., Todd, M. J., 1985. Identifying the set of always-active constraints in a system of linear inequalities by a single linear program.

Gelareh, S., 2008. Hub location models in public transport planning. Ph.D. thesis, Tecnical University of Kaiserslautern, Germany.

Gelareh, S., Maculan, N., Mahey, P., Monemi, R. N., 2013. Hub-and-spoke network design and fleet deployment for string planning of liner shipping. Applied Mathematical Modelling 37 (5), 3307–3321.

Gelareh, S., Nickel, S., 2011. Hub location problems in transportation networks. Transportation Research Part E 47, 1092–1111.

Gelareh, S., Nickel, S., Pisinger, D., 2010. Liner shipping hub network design in a competitive environment. Transportation Research Part E: Logistics and Transportation Review 46 (6), 991–1004.

Gendreau, M., Laporte, G., Semet, F., 1998. A branch-and-cut algorithm for the undirected selective traveling salesman problem. Networks 32 (4), 263–273.

Goldman, A., 1969. Optimal location for centers in a network. Transportation Science 3, 352–360.

Hakimi, S. L., 1964. Optimum locations of switching centers and the absolute centers and medians of a graph. Operations Research 12, 450–459.

Kara, B. Y., Taner, M. R., 2011. Hub location problems: The location of interacting facilities. In: Eiselt, H. A., Marianov, V. (Eds.), Foundations of location analysis. Springer, pp. 273–288.

Labbé, M., Rodríguez-Martin, I., Salazar-Gonzalez, J., 2004. A branch-and-cut algorithm for the plant-cycle location problem. Journal of the Operational Research Society 55 (5), 513–520.

Laporte, G., Nobert, Y., Arpin, D., 1986. An exact algorithm for solving a capacitated location-routing problem. Annals of Operations Research 6 (9), 291–310.

Maculan, N., Plateau, G., Lisser, A., 01 2003. Integer linear models with a polynomial number of variables and constraints for some classical combinatorial optimization problems. Pesquisa Operacional 23(1), 161 – 168.

Magnanti, T., Mireault, P., Wong, R., 1986. Tailoring Benders decomposition for uncapacitated network design. Netflow at Pisa, 112–154.

Magnanti, T., Wong, R., 1981a. Accelerating Benders decomposition: Algorithmic enhancement and model selection criteria. Operations Research 29 (3), 464–484.

Magnanti, T. L., Wong, R. T., 1981b. Accelerating benders decomposition: Algorithmic enhancement and model selection criteria. Operations Research 29(3), 464–484.

McDaniel, D., Devine, M., 1977. A modified Benders' partitioning algorithm for mixed integer programming. Management Science 24 (3), 312–319.

Nagy, G., Salhi, S., 1998. The many-to-many location-routing problem. Top 6 (2), 261–275.

Nickel, S., Schobel, A., Sonneborn, T., 2001. Hub location problems in urban traffic networks. In: Niittymaki, J., Pursula, M. (Eds.), Mathematics Methods and Optimization in Transportation Systems. Kluwer Academic Publishers, pp. 1–12.

O'Kelly, M., 1986a. The location of interacting hub facilities. Transportation Science 20 (2), 92–105.

O'Kelly, M. E., 1986b. Activity levels at hub facilities in interacting networks. Geographical Analysis 18 (4), 343–356.

O'Kelly, M. E., 1987. A quadratic integer program for the location of interacting hub facilities. European Journal of Operational Research 32 (3), 393–404.

O'Kelly, M. E., Bryan, D. L., 1998. Hub location with flow economies of scale. Transportation Research Part B: Methodological, Elsevier 32 (8), 605–616.

Poojari, C., Beasley, J., 2009. Improving benders decomposition using a genetic algorithm. European Journal of Operational Research 199 (1), 89–97.

Rodríguez-Martín, I., Salazar-González, J. J., Yaman, H., 2014. A branch-and-cut algorithm for the hub location and routing problem. Computers & OR 50, 161–174.

Sherali, H., Fraticelli, B., 2002. A modification of Benders' decomposition algorithm for discrete subproblems: An approach for stochastic programs with integer recourse. Journal of Global Optimization 22 (1), 319–342.

Toh, R. S., Higgins, R. G., 1985. The impact of hub and spoke network centralization and route monopoly on domestic airline profitability. Transportation Journal 24 (4), 16–27.

Wasner, M., Zäpfel, G., 2004. An integrated multi-depot hub-location vehicle routing model for network planning of parcel service. International Journal of Production Economics 90 (3), 403–419.

Yoon, M.-G., Current, J., 2008. The hub location and network design problem with fixed and variable arc costs: formulation and dual-based solution heuristic. J. Oper. Res. Soc. 59 (1), 80–89.